
\input amstex

\input epsf.tex

\documentstyle{amsppt}

\overfullrule=0pt

\define\sset{\subseteq}

\topmatter

\title\nofrills Canonical genus and the Whitehead doubles of pretzel knots\endtitle
\author  Mark Brittenham$ ^1$ and Jacqueline Jensen$ ^2$ \endauthor

\leftheadtext\nofrills{Mark Brittenham and Jacqueline Jensen}
\rightheadtext\nofrills{Canonical genus and doubles of pretzels}

\address   Department of Mathematics, University of Nebraska, Avery 203, Lincoln, NE 68588-0130  \endaddress
\email   mbrittenham2\@math.unl.edu \endemail

\address   Department of Mathematics and Statistics, Sam Houston State University, Huntsville, Texas 77341  \endaddress
\email jensen\@shsu.edu \endemail


\thanks  \hskip1pc$ ^{1}$ Research supported in part by NSF grant \# DMS$-$0306506 \endthanks
\thanks  \hskip1pc$ ^{2}$ Research supported in part by NSF grant \# DMS$-$0354281 \endthanks


\abstract 
We prove, for an alternating pretzel knot $K$, that the canonical
genus of its Whitehead doubles $W(K)$ is equal to the crossing number
$c(K)$ of $K$, verifying a conjecture of Tripp in the case of these
knots. 
\endabstract

\endtopmatter

\document

\heading{\S 0 \\ Introduction}\endheading

Every knot $K$ in the 3-sphere $S^3$ is the boundary of a compact
orientable surface $\Sigma\sset S^3$, known as a Seifert surface 
for the knot $K$. The first proof of this was given by Seifert [Se], who gave an 
algorithm which starts with a diagram of the knot $K$ and produces
such a surface. The algorithm consists of orienting the knot
diagram, breaking each crossing and reconnecting the resulting four ends
according to the orientation, without re-introducing a crossing, 
producing disjoint "Seifert" circles in the projection plane, bounding (after
offsetting nested circles) disjoint "Seifert" disks, and then reintroducing
the crossings  of $K$ by stitching the disks together with half-twisted
bands.

The minimum of 
the genera of the ``canonical'' surfaces built by Seifert's algorithm, over 
all diagrams of the knot $K$, is known as the {\it canonical genus} 
or {\it diagrammatic genus} of $K$, denoted $g_c(K)$. The minimum
genus over all Seifert surfaces, whether built by Seifert's
algorithm or not, is known as the {\it genus} of $K$, and
denoted $g(K)$.

Like all such invariants, defined as the minimum over an infinite 
collection of configurations, the canonical genus and the genus
can be extremely difficult to compute. But there are at least two
situations where such a determination can be reasonably carried out.
The first is when $g_c(K)=g(K)$, that is, a canonical surface has 
minimum genus among all Seifert surfaces. A candidate for such a surface
can, in principle (and usually in practice, as well), be verified to
have minimum genus using sutured manifold theory [Ga1]. This condition
cannot always be met, however; the genus need not equal the canonical genus,
for example with the Whitehead doubles of knots considered
here and elsewhere [Na],[Tr]. The second is
when Morton's inequality [Mo] is an equality. Morton's inequality
states that the $z$-degree of the HOMFLY polynomial $P_K(v,z)$ of
a knot $K$ is at most
twice the canonical genus $2g_c(K)$ of $K$. A
canonical surface $\Sigma$ whose genus is half of the $z$-degree,
if it exists, must therefore have genus equal to the canonical genus.
This condition also cannot always be met; the first examples 
where Morton's inequality was shown to be strict were
found by Stoimenow [St], and the authors [BJ] recently found several 
infinite families of examples.

Both methods succeed in computing the canonical genus of alternating
knots [Cr],[Ga2],[Mu], while the second
can compute the canonical genera of knots through 12 crossings [St]. 
The first method provides an approach to computing the canonical
genera of arborescent knots [Ga3]. In fact, the authors are aware of 
no example where
both approaches are known to fail. (The knots considered in [St].[BJ] can 
have their canonical genera computed by the first method, which is precisely why
the second method is known to fail.)

Recently, Morton's inequality has been applied to the computation
of the canonical genus of the Whitehead doubles $W(K)$ of knots $K$,
for which the first approach can shed no light; $g(W(K))=1$
for all non-trivial $K$. Tripp [Tr] computed the canonical genus of
the doubles of the $(2,n)$-torus knots $T_{2,n}$, showing that
$g_c(W(T_{2,n})) = n = c(T_{2,n})$ . For the proof he demonstrates 
(by induction on $n$) that the $z$-degree of the HOMFLY polynomial of the 
double is $2n = 2c(T_{2,n})$. A general construction, for any 
knot $K$, provides a canonical surface with genus $c(K)$ (we 
reproduce this construction in the next section). So we have the inequalities

\smallskip

\centerline{$2c(T_{2,n}) = 2n = $deg$_z P_{W(T_{2,n})}(v,z) 
\leq 2g_c(W(T_{2,n}))\leq 2c(T_{2,n})$}

\smallskip

\noindent from which $g_c(W(T_{2,n})) =c(T_{2,n})$ follows.
Nakamura [Na] has extended this argument to prove that for 2-bridge knots
$K$, $g_c(W(K)) = c(K)$. The main part of the argument is again an
inductive proof that deg$_z P_{W(K)}(v,z) = 2c(K)$ .


\leavevmode

\epsfxsize=2.5in
\centerline{{\epsfbox{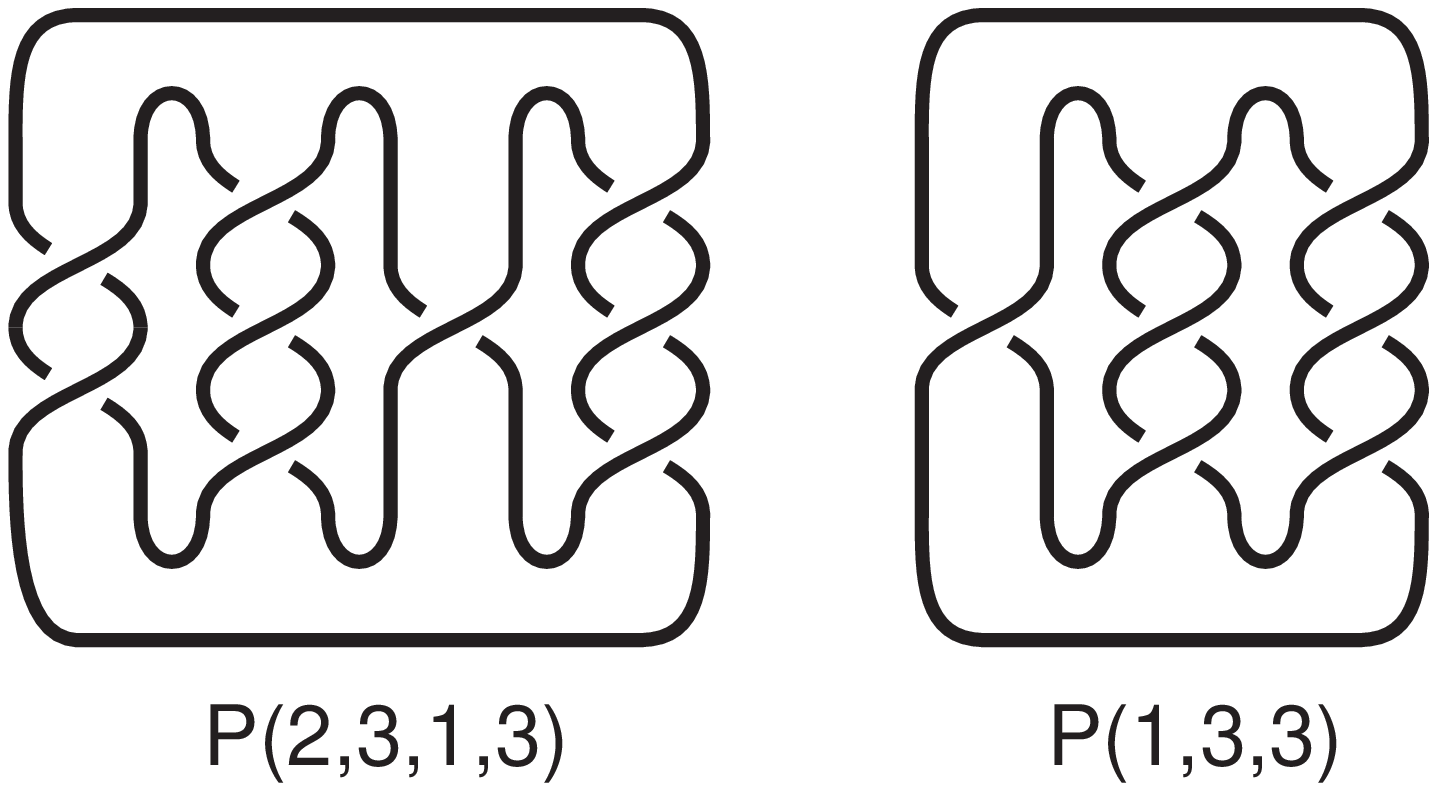}}}

\centerline{Figure 1: Pretzel knots}

\medskip

In this paper we provide another extension of Tripp's result, to 
alternating pretzel knots $P(k_1,\ldots k_n)$, $k_1,\ldots ,k_n\geq 1$.
Note that $P(k_1,\ldots k_n)$ is a knot iff either $n$ is odd and at 
most one $k_i$ is even, or $n$ is even and exactly one $k_i$ is even.

\proclaim{Theorem 1}
If $K$ is a pretzel knot $P(k_1,\ldots k_n)$ with $k_1,\ldots ,k_n\geq 1$,
then $g_c(W(K)) = k_1+\cdots +k_n = c(K)$ .
\endproclaim

The main tool in the proof is a proposition which 
gives a method for building new knots $K$ satisfying 
$2c(K) = $deg$_z P_{W(K)}(v,z)$ from old ones $K^\prime$.

\proclaim{Proposition 2}
If $K^\prime$ is a knot satisfying deg$_z P_{W(K^\prime)}(v,z) = 2c(K^\prime)$,
and if for a $c(K^\prime)$-minimizing diagram for $K^\prime$ we replace a
crossing, thought of as a half-twist, with three half-twists
(see Figure 2), producing a knot $K$, then 
deg$_z P_{W(K)}(v,z) = 2c(K)$, and 
therefore $g_c(W(K)) = c(K)$ .
\endproclaim

\medskip

Theorem 1 follows immediately by repeated application of Proposition 2 to the
crossings of the $(2,n)$ torus knots $=P(1,\ldots ,1)$ 
(when all $k_i$ are odd) and
the 2-bridge knots $[2,n-1]=P(2,1,\ldots ,1)$ (when one of the $k_i$ 
(wolog $k_1$) is even), where the initial hypothesis is satisfied by the 
results of Tripp and Nakamura. It should be clear that Theorem 1 is not the strongest
result that can be proved in this way; the pretzel knots are just the simplest
new class of knots encountered. A more general statement is:


\leavevmode

\epsfxsize=2.5in
\centerline{{\epsfbox{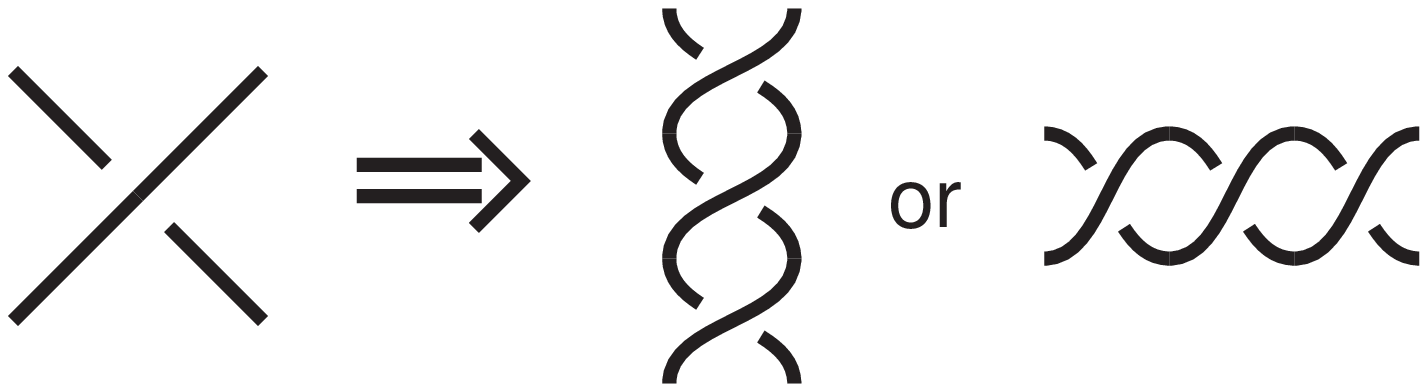}}}

\centerline{Figure 2: Introducing twists}

\medskip

\proclaim{Proposition 3} 
Let $\Cal K$ be the class of knots having diagrams which can 
be obtained from the standard diagram of the left- or right-handed trefoil knot $T_{2,3}$ by repeatedly 
replacing a crossing, thought of as a half-twist, by a full twist. Then
for every $K\in {\Cal K}$, deg$_z P_{W(K)}(v,z) = 2c(K)$, 
and so $g_c(W(K)) = c(K)$ .
\endproclaim

Using Tripp's and Nakamura's results, we could have stated this corollary using the
a priori larger class obtained by starting with the reduced alternating
diagrams of any 2-bridge knot. But we will show that all 2-bridge knots already
lie in $\Cal K$, and so this ``larger'' class would really just be $\Cal K$. Stated
differently, Proposition 3 implies the results of Tripp and Nakamura!
We should note that the knots built from the left-handed trefoil are the mirror images
of those built from the right-handed trefoil, and the $z$-degree of the HOMFLY polynomial
is unchanged under taking mirror images. So it suffices to prove the result for only
one of the two collections of knots.

\heading{\S 1 \\ Notations and preliminaries}\endheading

$K$ will always denote a knot or link in the 3-sphere $S^3$, 
$N(K)$ a tubular neighborhood of $K$, $E(K) = S^3 \setminus$ int$N(K)$ the
exterior of $K$, and $\Sigma$ a Seifert surface for $K$, which we treat as embedded
in $S^3$. In this paper we take a diagrammatic approach to the construction of 
Whitehead doubles $W(K)$ of a knot $K$. Given a diagram $D$ of a knot $K$, 
we construct diagrams for the Whitehead doubles of $K$ as follows:


\leavevmode

\epsfxsize=2.5in
\centerline{{\epsfbox{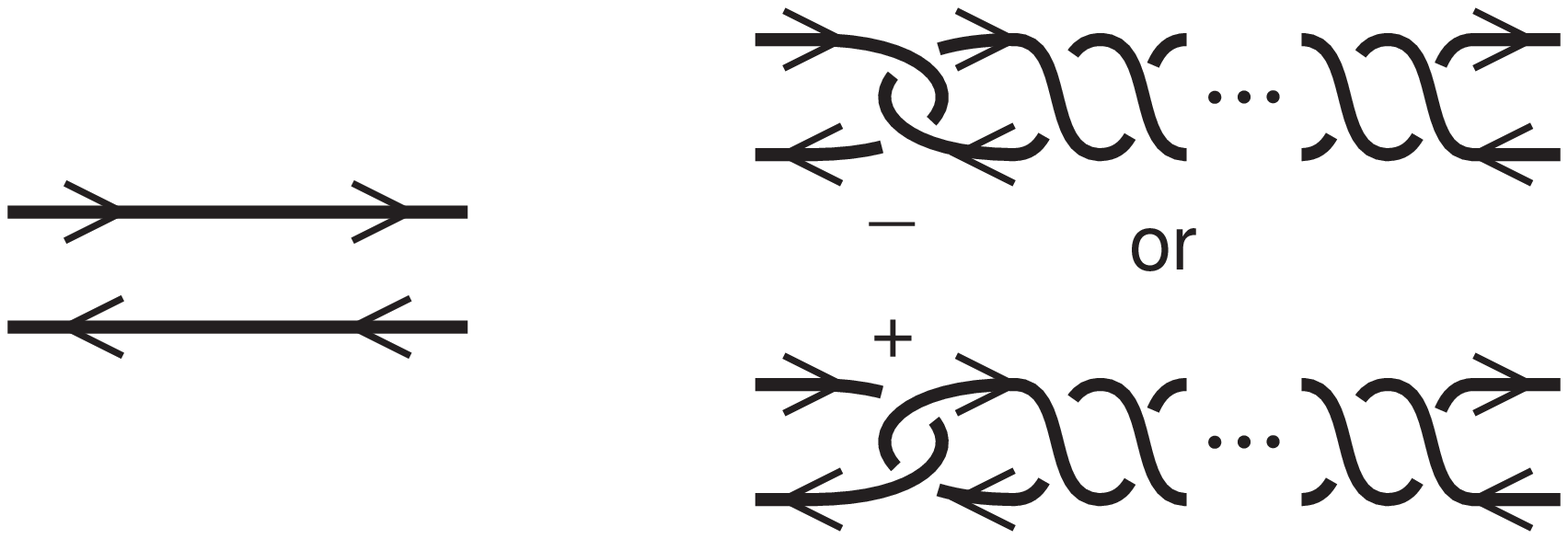}}}

\centerline{Figure 3: Whitehead doubles}

\medskip

Starting with the diagram $D$, choose an orientation for the underlying knot. Draw
a second, parallel, copy of $K$ pushed off of $K$ to the right, according to the
orientation of $K$, and oriented in the opposite direction, giving an oriented link $D(K)$ 
(which for lack of a better term, we will
call the {\it flat double} of $K$). Note that we will always treat $D(K)$ as an
{\it oriented} link, with the orientation just described. 
$W(K)$ is obtained from $D(K)$ by replacing
a pair of parallel arcs in the diagram by a ``twisted clasp'', as in Figure 3.
These twisted clasps consist of either a right-handed (+) or left-handed (-)
clasp together with $|n|$ full twists to the parallel arcs (right-handed for
$n>0$, left-handed for $n<0$). If we need to specify a Whitehead double, rather
than deal with them as a class, we will denote it by $W(K,n,+)$ or $W(K,n,-)$ .
The interested reader may satisfy him/herself that this description of Whitehead
doubles is equivalent to any other description he or she has encountered. Note
that the notation $W(K,n,\pm)$ is very diagram-dependent; a different diagram
will build the same collection of knots, but might label them differently
(depending upon the writhe of the diagram).

For a given non-trivial knot $K$, if we build its Whitehead doubles using a $c(K)$-minimizing
diagram as above, then after ``hiding'', by isotopy, the full twists of $W(K)$ inside of 
one of the small squares formed from a crossing of $K$, as in Figure 4a, Seifert's algorithm
will build a canonical surface for $W(K)$ of genus $c(K)$. To see this, note that the canonical
surface is a checkerboard surface; this was the point to building our Whitehead doubles with 
full twists in the parallel strands, and then hiding them inside the doubling of a crossing of 
$K$. The local pictures are as in Figure 4b.

\leavevmode

\epsfxsize=4in
\centerline{{\epsfbox{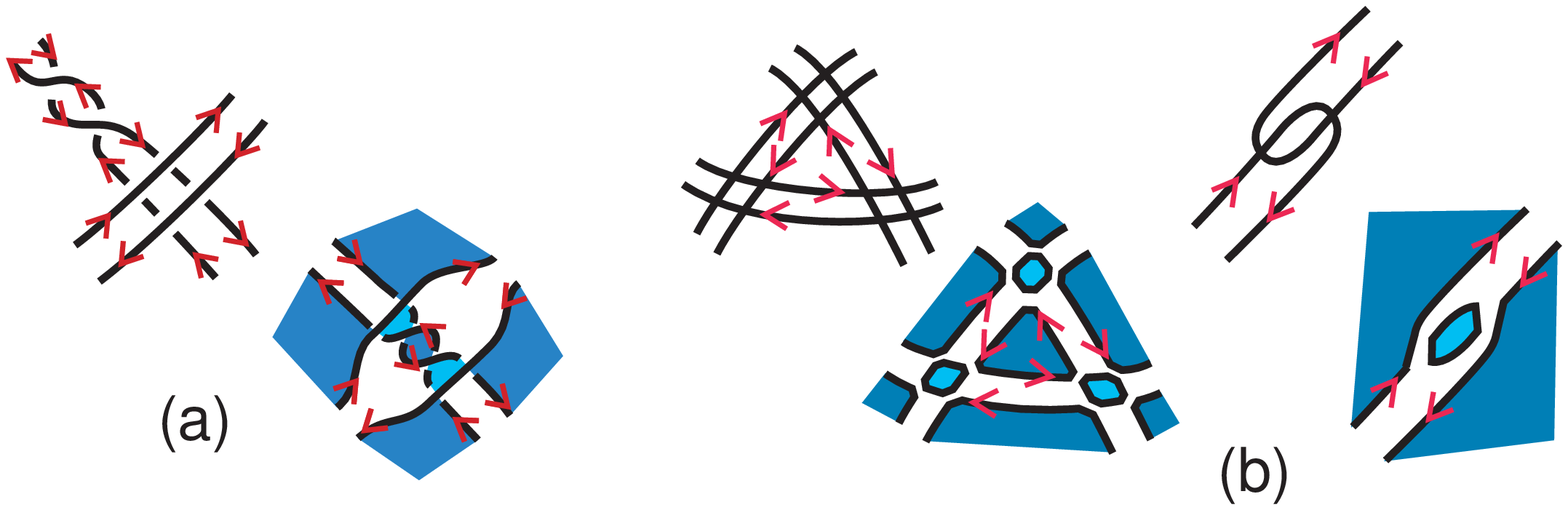}}}

\centerline{Figure 4: A canonical surface for the Whitehead double}

\medskip

Continuing our calculation of the genus, note that the canonical surfaces for 
all of the Whitehead doubles are homeomorphic,
and so have the same genus; the surfaces differ only in the direction of twisting
at the clasp and the number of full twists to the Seifert disk we ``hid'' the 
full twists of the doubling in. So we focus on the surface for the Whitehead double
with no full twists, $W(K,0,+)$ . Its canonical Seifert surface $\Sigma$ has one Seifert disk 
for each complementary region of the diagram $D$ of $K$, one for each crossing of $K$,
and one for the Whitehead doubling clasp. It has one half-twisted band for each crossing of
$W(K,0,+)$, that is, $4c(K)+2$ twisted bands. So 
$1-2g(\Sigma) = \chi(\Sigma) = |$Seifert disks$|$ - $|$bands$|$ 
= ($c(K)$ + $|$regions of $D|+1$) - ($4c(K)+2$) = $|$regions$| - 3c(K) -1$. But since 
$D\subseteq S^2$ is a 4-valent graph whose complementary regions are all disks, we have
(using the usual notation) $v-e+f=2$, and $2e=4v$ (since all vertices have valence $4$),
so $e=2v$, and so $f$ = $|$regions$|$ = $2+v = 2+c(K)$. So 
$1-2g(\Sigma) = 2+c(K) -3c(K) -1 = 1-2c(K)$ . So $g(\Sigma) = c(K)$, as desired.

\medskip

This implies that $g_c(W(K))\leq c(K)$ for all nontrivial knots $K$. Using 
Morton's inequality, we can establish the reverse inequality if we can 
show that deg$_z P_{W(K)}(v,z) = 2c(K)$ . Our goal is to prove this equality
for a class of knots that includes the alternating pretzel knots.
Central to our proof is therefore a computation of the $z$-degree of the HOMFLY polynomial.

\smallskip

The HOMFLY polynomial [FHLMOY] is a 2-variable Laurent polynomial defined 
for any oriented link, and may be thought of as the unique polynomial
$P_K(v,z)$,
defined on link diagrams and invariant under the Reidemeister moves,
satisfying $P_{\text unknot}(v,z)=1$, and $v^{-1}P_{K_+} - vP_{K_-} = zP_{K_0}$,
where $K_+,K_-,K_0$ are diagrams which all agree except at one crossing, where
they are given as in Figure 5. Here we are following Morton's convention in
the naming of variables. 

\leavevmode

\epsfxsize=2in
\centerline{{\epsfbox{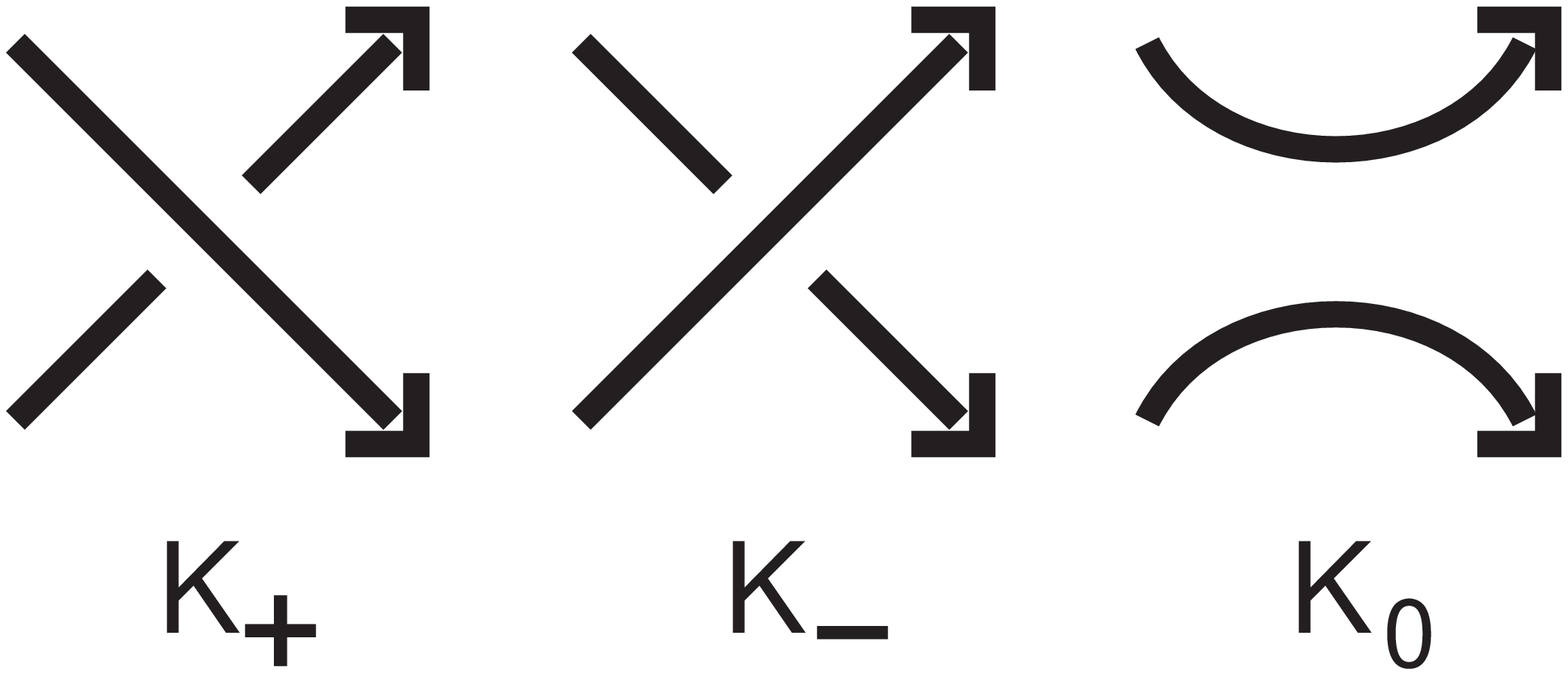}}}

\centerline{Figure 5: HOMFLY polynomial}

\medskip

This skein relation gives an inductive method for computing the HOMFLY polynomial,
since one of $K_+,K_-$ will be ``closer'' than the other to the unlink,
in terms of unknotting number, while $K_0$ has fewer crossings.
It also allows any one of the polynomials to be computed from the other two.

If we write $P_K(v,z)$ as a polynomial in $z$, having coefficients which are
polynomials in $v$, then Morton showed [Mo] that, 
for any oriented diagram $D$ of an oriented knot or link $K$,
the $z$-degree $M(K)$ of
$P_K(v,z)$ satisfies $M(K)\leq  c(D)-s(D)+1$, where $c(D)$ is the number of crossings of the
diagram and $s(D)$ is the number of Seifert circles of $D$. Since the
Seifert surface $\Sigma$ built by Seifert's algorithm from this diagram has
Euler characteristic $s(D)-c(D) = \chi(\Sigma) = (2-|K|)-2g(\Sigma)$,
where $|K|$ is the number of components of $K$, we have 
$M(K)\leq c(D)-s(D)+1 = 1-\chi(\Sigma) = 2g(\Sigma) + (|K|-1)$ for every
canonical Seifert surface for $K$. Consequently, for a knot $K$, 
$M(K)$ is bounded from above by twice the canonical genus of $K$, $2g_c(K)$. 
Our main interest in the HOMFLY polynomial will therefore be in a computation 
of this $z$-degree.
In particular, since the degree of the sum of two polynomials cannot exceed the 
larger of their two degrees, and is equal to the larger if the degrees are unequal, 
we get the basic inequalities (letting $K_+$ denote
the knot with diagram $D_+$, etc.) for the $z$-degree of $P_K$, $M(K)$ :

\smallskip

\hskip.2in $M(K_+)\leq \max\{M(K_-),M(K_0)+1\}$

\hskip.2in $M(K_-)\leq \max\{M(K_+),M(K_0)+1\}$

\hskip.2in $M(K_0)\leq \max\{M(K_+),M(K_-)\}-1$

\smallskip

\noindent where equality holds if the two terms in the maximum are unequal, from the
equations $P_{K_+}=v^2P_{K_-}+vzP_{K_0}$, $P_{K_-}=v^{-2}P_{K_+}-v^{-1}zP_{K_0}$, and
$P_{K_0}=v^{-1}z^{-1}P_{K_+}-vz^{-1}P_{K_-}$. For example, since if we change one of the crossings
in the clasp of $W(K) = L_\pm$ we get the unknot $= L_\mp$, while if we break the crossing we get the 
link $D(K) = L_0$, we have $M(D(K)) \leq \max\{M(W(K)),0\} -1$, with equality if
$M(W(K))>0$. In particular, if we already know that $M(D(K)) > 0$, then $M(D(K))=M(W(K))-1$.

\heading{\S 2 \\ Proof of the main result}\endheading

In this section we prove Proposition 4 below, from which Propositions 2 and 3 and Theorem 1 will follow. 
Our main goal is to compute the $z$-degrees of the HOMFLY polynomials of Whitehead 
doubles $W(K,n,\pm)$ of knots $K$.
First, we get rid of the clasp. As we just saw, $M(D(K) )\leq \max\{M(W(K)),0\} -1 = M(W(K))-1$.
So to prove Proposition 3 we wish to establish that if a knot $K$ has a diagram $D$
obtained from the standard diagram $D_0$ of the trefoil
knot $T_{2,3}$ by repeatedly replacing crossings, treated as half-twists, by full twists,
then for a doubled link with $n$ full twists $D(K,n)$, we have $M(D(K,n))=2c(K)-1$. 
We note that, by induction, this replacement process always builds
alternating diagrams with no nugatory crossings. A nugatory crossing is a crossing
in the diagram $D$ so that some pair of opposite complementary regions represent the same connected
component of the complement of the projection of $K$. Such a phenomenon is preserved when passing 
between a diagram with a full twist and one with a half-twist; 
if the crossing is one of those in the full twist, then the half-twist will be
nugatory, while if not, then the nugatory crossing disjoint from the half- or full
twist will persist.
So, in our arguments, we will always have $c(K)=c(D)$, i.e., $D$
is a diagram with minimal crossing number. So what we will show is that $M(D(K))=2c(D)-1$.

Next we get rid of the $n$ full twists. Given a specific double $D(K,n)$ of knot $K$,
changing one of the crossings among the full twists yields $D(K,\pm(|n|-1))$ (after isotopy),
while breaking the crossing yields the unknot $L$, with $M(L)=0$.
So our basic inequalities yield 
$M(D(K,\pm(|n|-1)))\leq \max\{M(D(K,n)),1\}$, with equality so long as
$M(D(K,\pm(|n|-1)))\neq 1$. By induction, then, since for a non-trivial knot
$K$ we have $c(K)\geq 3$, so long as
we show that $M(D(K,0))=2c(D)-1$ we will have established that $M(D(K,n))=2c(D)-1$
for all $n$. Since we will focus on the case $n=0$, from this point we will use
$D(K)$ to denote the flat double $D(K,0)$.

In general, using the same argument as in the last section, if $D^\prime$ is a diagram 
of a link $L$ with doubled
link $D(L)$, since the projection of $L$ is a graph with $c(D^\prime)$ vertices and 
$c(D^\prime)+2$ faces, 
Morton's inequality gives 

\smallskip

\noindent $M(D(L))\leq c(D(D^\prime))-s(D(D^\prime))+1$ 

\hfill $= (4c(D^\prime))-(c(D^\prime)+(c(D^\prime)+2))+1=2c(D^\prime)-1$.

\smallskip

\noindent This basic estimate will be used several times in the proofs that follow.

\medskip

Showing that $M(D(K))=2c(D)-1$ has the advantage of making sense when the underlying 
$K$ is a link, and not just a knot; the Whitehead doubling process does not really apply to a link. 
It is therefore a better setting to use for an inductive proof;
replacing a crossing in a knot diagram with a full twist produces a 2-component
link, not a knot.
We will therefore work in the setting of flat doubles $D(K)$ of knots or links $K$, following
the orientation conventions established in section 1; the parallel doubled strands
for every component of $K$ are oriented oppositely, so that the parallel strand is to the right.
What we show is:

\proclaim{Proposition 4}
If $L$ is a non-split link with a diagram $D^\prime$ satisfying  $c(D^\prime)=c(L)$ 
and $M(D(L))=2c(D^\prime)-1$, and $K$
is a link having diagram $D$ obtained by replacing a crossing in the diagram $D^\prime$ 
with a full twist (so that $c(D)=c(D^\prime)+1$), then $M(D(K))=2c(D)-1$.
\endproclaim


From the discussion at the beginning of this section describing how
information about $M(D(K))$ translates to information about $M(W(K))$, 
Proposition 4 implies Proposition 3
by induction on the number of full twists introduced, once we establish the base case.
But a direct calculation (see also [Tr],[Na]) establishes that
$M(D(T_{2,3})) = 5 = 2\cdot 3-1 = 2c(T_{2,3})-1$, giving the base case for the induction.
Since at every stage the process of introducing full twists builds an 
alternating diagram $D$ for $K$ with no nugatory crossings, $c(K)=c(D)$ at every stage of
the induction, and since the link projection is always connected, Menasco [Me]
implies that the underlying link is non-split.
So Proposition 4 applies at every step of the induction.
Proposition 2 follows by applying Proposition 4 twice. 
And as remarked in the introduction,
Theorem 1 follows from Proposition 2, using the two base cases the torus 
knots $T_{2,n}=K(1,\ldots ,1)$, supplied by [Tr], and the twist knots $K(2,1,\ldots ,1)$, 
supplied by [Na]. We should note that the pretzel knots $K(n_1,n_2)$ are the torus knots
$T_{2,n_1+n_2}$.

We now turn to the proof of Proposition 4. The proof is 
by induction on $c(K)$, and is essentially based on a large skein tree
diagram calculation, with the diagram $D(D)$ for $D(K)$ at the root and $D(D^\prime)$ as one of the leaves.
We will show that under the hypotheses of the proposition we have $M(D(K))= M(D(L))+2$, which
establishes the inductive step. Since we are only interested in 
$M(D(K))=$maxdeg$_z(P_{D(K)})$, we will only track this quantity through the calculation, rather than the 
entire HOMFLY polynomial, using the inequalities established at the end of section 1. 

\leavevmode

\epsfxsize=5in
\centerline{{\epsfbox{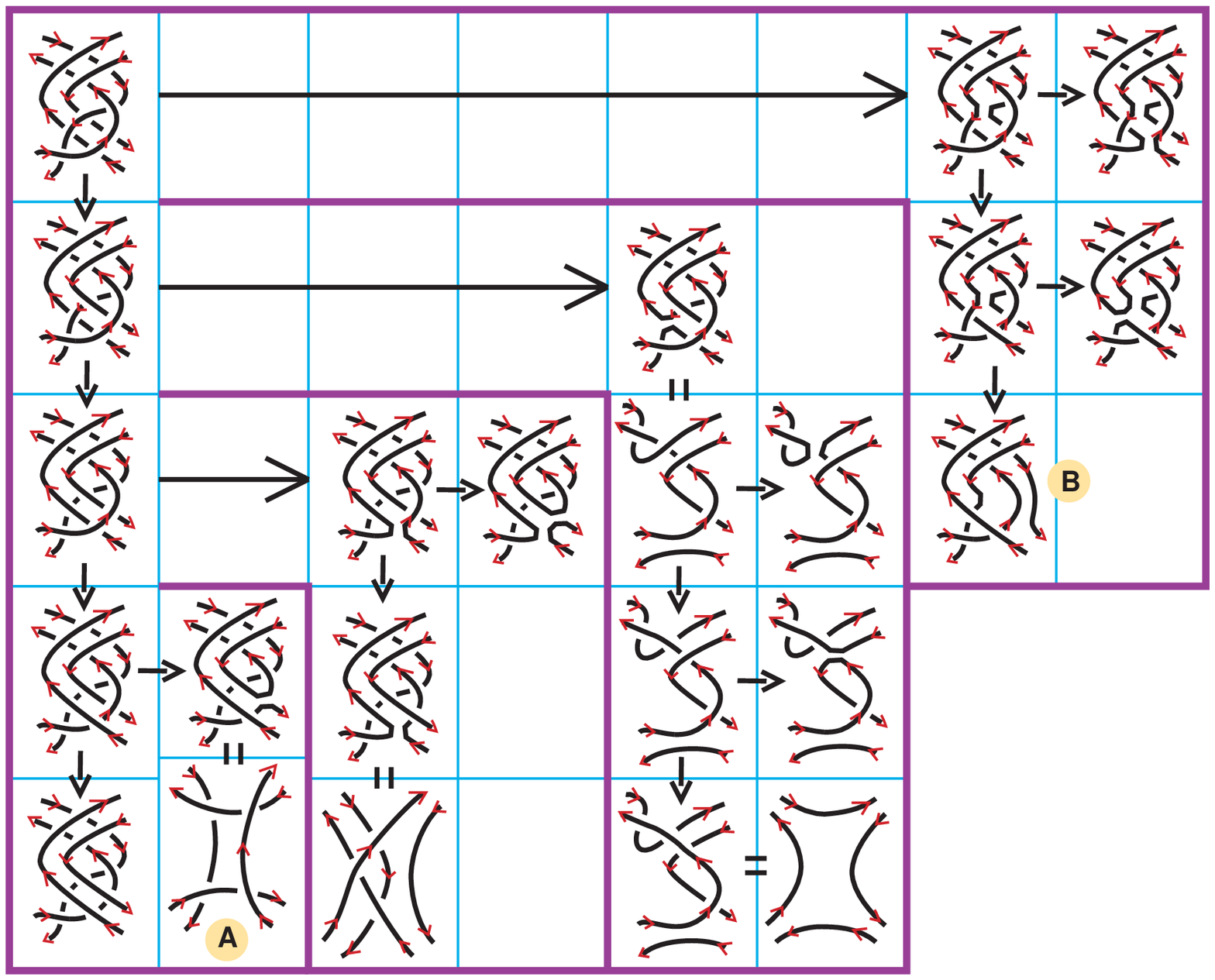}}}

\centerline{Figure 6: Main skein tree diagram}

\medskip

\medskip

The main skein tree diagram that we employ is given in Figure 6. The upper left corner of
the figure is a portion of our diagram $D(D)$ for $D(K)$ at the full twist of $K$ which
exists by hypothesis, with an additional pair of canceling half-twists
introduced on either side of a quartet of overcrossings. A horizontal step in the 
diagram involves breaking a crossing, according 
to the orientation of the strands; a downward step involves changing the same crossing. 
Except for the addition
of a full-twist to the two parallel strands, the upper right corner of the figure is the corresponding
portion of the diagram $D(D^\prime)$ of $D(L)$. The two right arrows across the top are the basis for
our induction; they imply, essentially, that $P_{D(K)} = z^2v^{-2}P_{D(L)} +$ other terms. The proof involves
showing that the other terms do not have high enough $z$-degree to interfere with the
conclusion that $M(D(K))=M(D(L))+2$. This demonstration will occupy the remainder of the
section. 

This in turn involves arguing up from the bottom of the 
skein tree, using the inequalities at the end of section 1 to bound the $z$-degrees of
these other terms from above. The initial steps in this process are done using Morton's inequality, 
since we can compute the genus of the Seifert surface built from the diagram given to
us at the leaves of the tree, or, more simply, the difference between this and the genus of 
the surface built from $D(D)$,
by computing the changes in the number of crossings (which is evident) and the number of
Seifert circles (which we can do since we know precisely where these circles come from for $D(D)$).
In the generic case, this will be precisely how our estimates are made. But there is also the possibility
that at some of the leaves our link will include a (possibly) twisted double of the unknot as a split
component, and our estimates will then be off by the $z$-degree of its HOMFLY polynomial; they will be too
low by 2 for each such component. The estimates need to be refined slightly in this case; for ease of exposition,
we deal with this possibility at the end of the proof, and, in what follows now, act as 
if such a split component has not been produced.

To start, we make note that the local picture in the upper left corner of the skein diagram,
if we were to undo the canceling pair of half-twists (since it is this simpler diagram that we 
are postulating minimizes the quantity $c(D(D))-s(D(D))+1$),
contains 8 crossings of $D(D)$ and contains parts of 7 Seifert disks - 2 from the crossings of $K$ and
5 from the complementary regions of $K$. Note that these 5 regions of $K$ are distinct; the one in the middle
of the full twist is demonstrably so from the figure, and the other four come from the four regions surrounding a
single crossing of $K^\prime$, which are all distinct because this crossing is not nugatory
(that is, opposite regions are distinct); otherwise, $c(D^\prime)$ would be greater 
than $c(L)$. As remarked earlier, replacing a crossing by a full twist cannot introduce a 
nugatory crossing. It is therefore this contribution, 8 crossings and 7 disks, that we will compare
with the leaves of our skein tree to compute what genus surface each diagram will build, giving an
upper bound, by Morton's inequality, on the $z$-degree of the HOMFLY polynomial 
of the underlying oriented link.

\leavevmode

\epsfxsize=3in
\centerline{{\epsfbox{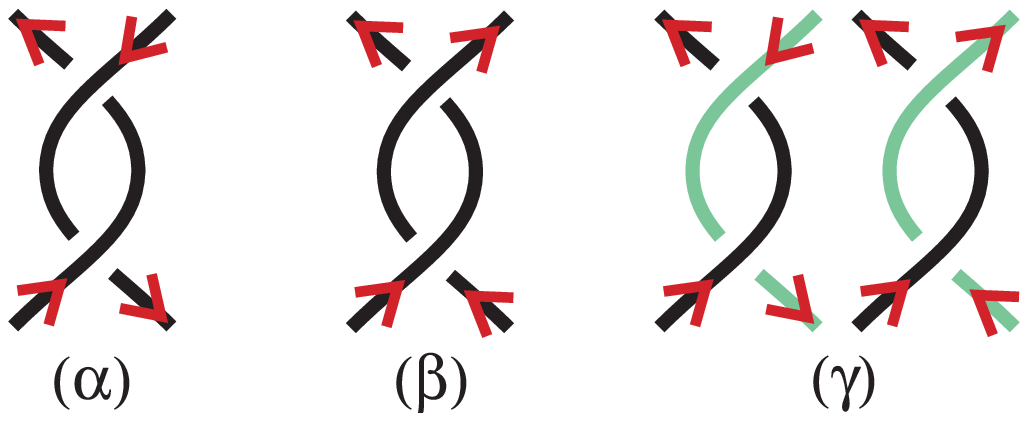}}}

\centerline{Figure 7: Three cases}

\medskip

Our argument will really be three parallel arguments, based on the underlying orientations,
and number of components, of the knot or link $K$ involved at the full twist where we
are carrying out our skein tree argument. There are essentially
three cases, given in Figure 7; ($\alpha$) and ($\beta$) are when the two strands at the full twist belong 
to the same component of $K$, with the opposite, respectively parallel, orientation, while ($\gamma$) is
when the two strands lie in different components.
We must at several points break our argument into three parallel ones,
according to which of these three cases we find ourselves in, when, 
 the diagrams as given at the leaves of the tree will typically yield canonical
surfaces whose genera are too high for our purposes. In each case we will wish to get rid of some 
of the crossings 
(thought of as half-twists) 
from the local diagrams at the leaves of the tree, in order to simplify the diagram and give us
a lower genus Seifert surface. So we will push each half-twist through the rest of the
link diagram (which consists entirely of pairs of parallel strands, having been untouched by the
skein moves that we have carried out), along these parallel strands, 
which will return the crossing to another one of the four corners of 
the local diagram. This is precisely where the three cases arise; each behaves differently
under this operation. 

In case ($\alpha$), pushing a crossing
from the upper left corner will return it to the upper right (and vice versa), while from the lower
left it will return to the lower right. In case ($\beta$), pushing from the upper left will
return to the lower right, and from the upper right will return to the lower left. In case
($\gamma$), pushing from the upper left returns to the lower left, and from the upper right
returns to the lower right. The other point to note is that if we push from two corners in
such a way that they return to the \underbar{other} two corners, then the net effect outside of the 
local diagram will be that all strands appear to be the same, but all of their orientations
will be reversed. (If there are entire pairs of parallel components of $D(K)$ outside of the local diagram,
we also interchange them by an isotopy, so that this is true for \underbar{all} components outside of 
the local diagram.) The point to this is that as a result, all of the Seifert circles outside
of the local diagram, when running Seifert's algorithm, are the same after this pair of crossing 
pushes; they have simply had all of their orientations reversed. In particular, in our count of the
number of Seifert circles at the leaves of the skein tree, we may act, after these pushes, as
if the Seifert disks outside still come precisely from the crossings of $K$ and the complementary
regions of $K$, as before in Figure 4 above.

\leavevmode

\epsfxsize=4.5in
\centerline{{\epsfbox{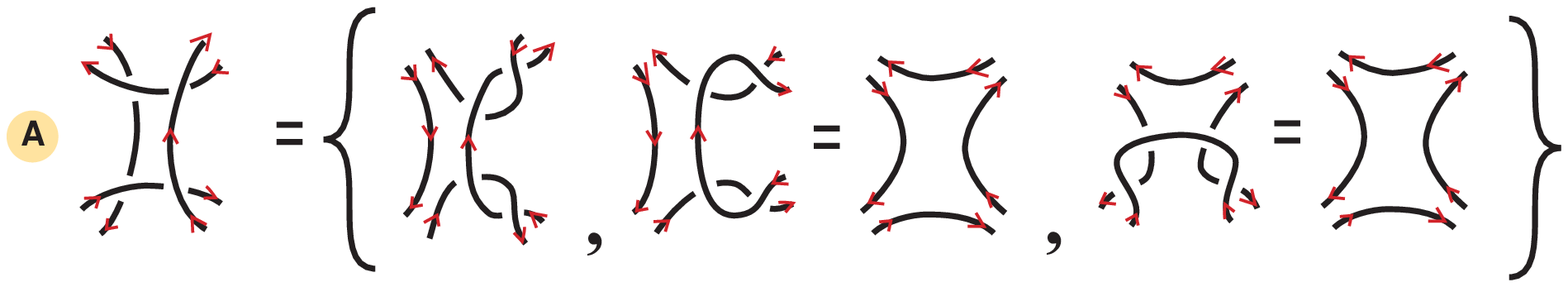}}}

\centerline{Figure 8: The ``A'' cases}

\medskip

In the case of the spot marked ``A'' in Figure 6, we have four extra crossings at the corners which if
we leave in place will give a surface with genus which is too high for our purposes when we run Seifert's
algorithm. Figure 8 shows the result of pushing a pair of the crossings (upper and lower left in the
first two cases, upper right and left in the third). The results are identical, except that in the first 
case we have two extra pairs of full twists. We will show shortly that these full twists, like our argument
at the beginning of this section passing from $D(K,n)$ to $D(K,0)$, do not affect the $z$-degree, but first let us analyze the $z$-degree
in the second and third cases. These diagrams $D_A$ have 8 fewer vertices than the local diagram for $D(D)$, and 3 fewer Seifert
disks. Therefore, for these diagrams, the resulting Seifert surface $\Sigma^\prime$ has
$c(D_A)=c(D(D))-8$ and $s(D_A)=s(D(D))-3$. By Morton's inequality, this yields
$M(D_A)\leq c(D_A)-s(D_A)+1=(c(D(D))-s(D(D))+1)-5=(2c(K)-1)-5=2c(K)-6$ for the last two cases. 
This, we shall see, is a $z$-degree that is too low to interfere with our main calculation.

Removing the pair of full twists in the first case
will give a diagram identical to the other two. The effect on the $z$-degree that these twists have
can be determined from the skein tree in the left of Figure 9; each of the two times we break the crossing,
the resulting short arc can be pushed out of the local diagram to return, since we are in case ($\alpha$),
to the other top or bottom corner. When it does, we see that the resulting link is another double
$D(L)$ of some link $L$ (possibly with a full twist). The full twist can be removed by the method of
the beginning of this section, without changing the $z$-degree, unless $L$ is the unknot, in which case we will
have degree 2. Since our original knot or link $K$ had no nugatory
crossings, however, the isotopy of the short arc has erased at least two crossings, in addition to 
the two crossings in our local diagram which have now disappeared (Figure 9, middle); the arc of
$K$ between the two corners must meet other crossings,
else it represents a monogon in the digram, and we have a nugatory crossing. 

\leavevmode

\epsfxsize=4.5in
\centerline{{\epsfbox{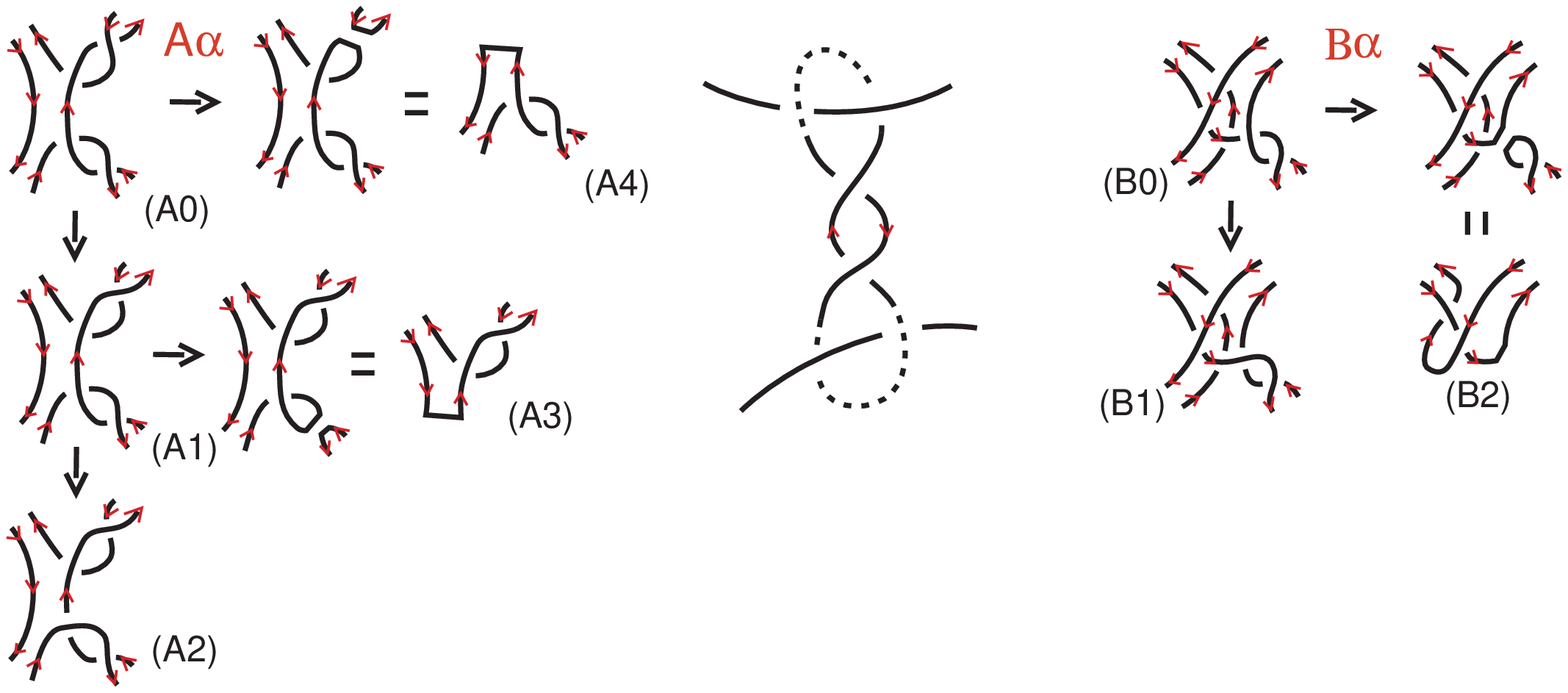}}}

\centerline{Figure 9: Removing full twists}

\medskip

So $L$ has at
least four fewer crossings than $K$, so $M(D(L)) \leq 2c(L)-1 \leq 2c(K)-9$. Consequently, 
as we move up the skein tree on the left of Figure 9, we have, since the link $K_{A2}$ is the same link 
found in the second and third cases, 
$M(K_{A1})\leq \max\{M(D_{A2}),M(D_{A3})+1\}\leq\max\{2c(K)-6,2c(K)-8\}=2c(K)-6$,
so $M(K_{A0})\leq \max\{M(K_{A1}),M(D_{A4})+1\}\leq\max\{2c(K)-6,2c(K)-8\}=2c(K)-6$.
So in all three cases we have $M(D_A)\leq 2c(K)-6$.

\leavevmode

\epsfxsize=4.5in
\centerline{{\epsfbox{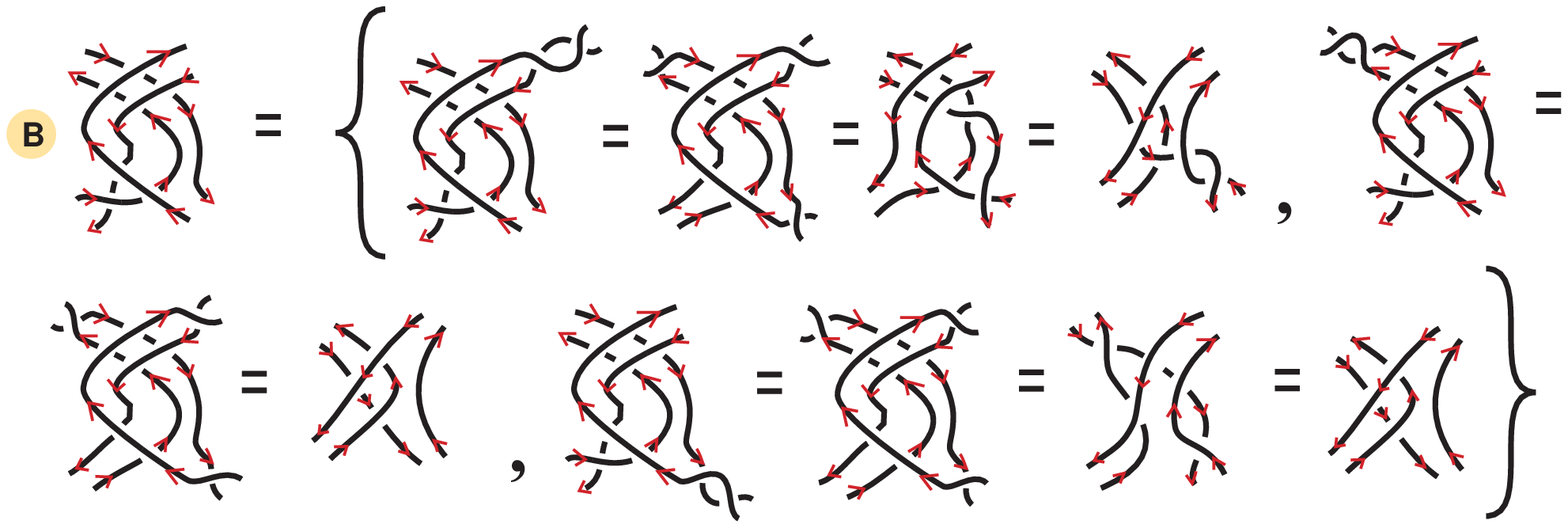}}}

\centerline{Figure 10: The ``B'' cases}

\medskip

In a similar vein, we analyze the $z$-degree for the link with diagram at spot ``B'' in our skein tree. The situation for the
three cases is illustrated in Figure 10. In this situation we need to ``invent'' some half-twists to push, since the diagram 
does not obviously have a pair of half-twists to work with. When we do so, we obtain the diagrams in the figure. Again, 
the first case has an extra full twist that we must deal with, but first we look at the last two cases. Their local
diagrams are identical and, compared with $D(D)$, $D_B$ has 5 fewer crossings and 2 fewer Seifert disks. So, by
Morton's inequality, 
$M(D_B)\leq c(D_B)-s(D_B)+1=((c(D(D))-5)-(s(D(D))-2)+1)=(c(D(D))-s(D(D))+1)-3=(2c(K)-1)-3=2c(K)-4$ for the last two cases.

In the first case we remove the full twist in the same manner that we did for spot ``A''; this is illustrated on the
right in Figure 9. Again we find that breaking a crossing of the full twist yields a double of a link with at least
four fewer crossings (with a full twist), and so $M(K_{B2})\leq 2c(K)-9$. The link given by the diagram $D_{B1}$
is not exactly the same as in the last two B-cases - one of the crossings is different - but it has the same
projection to the plane and so has the same number of crossings and the same Seifert disks. So $M(K_{B1})\leq 2c(K)-4$.
Working our way up the skein tree then yields $M(K_{B0})\leq \max\{M(K_{B1}),M(K_{B2})+1\}\leq
\max\{2c(K)-4,2c(K)-8\}=2c(K)-4$. So in all three cases we have $M(D_B)\leq 2c(K)-4$.

\leavevmode

\epsfxsize=5in
\centerline{{\epsfbox{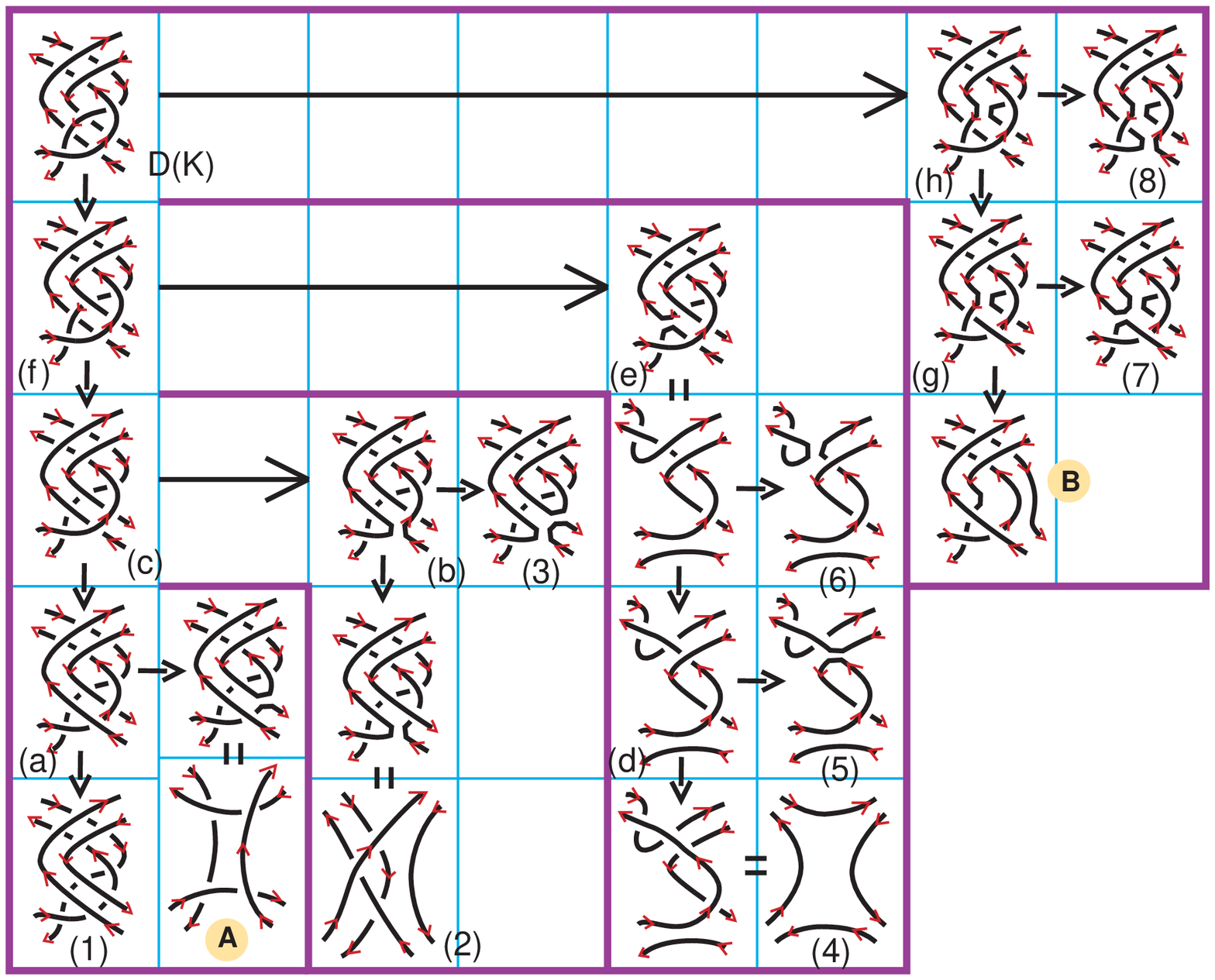}}}

\centerline{Figure 11: Main skein tree diagram, labeled}

\medskip

To proceed with our degree calculation, we need similar estimates for all of the other leaves of 
our skein tree diagram. We reproduce our main skein tree diagram as Figure 11, with additional labels
to mark the nodes of the tree. The link at
the spot marked $(1)$ is isotopic to the double of a link diagram $D^\prime$ with
two fewer crossings than $K$, so $M(K_1)\leq 2c(D^\prime)-1=2c(K)-5$. At $(2)$ we have the diagram
encountered in the (first case of) spot ``B'', so $M(K_2)\leq 2c(K)-4$. At $(4)$, we have the 
diagram encountered in case ``A'',
so $M(K_4)\leq 2c(K)-6$. At $(5)$, we have, after isotopy, the double of a link diagram with
2 fewer crossings than $K$, so $M(K_2)\leq 2c(D_5)-1=2c(K)-5$. 

At $(3)$, $(6)$ and $(7)$, we have diagrams
with short arcs in one corner that we can push out of the local diagram to return at one of the other three
corners, after which, as in the arguments for spots ``A'' and ``B'', we have the double of a link
with at least 2 fewer crossings; the crossings that were in the local diagram are now gone. (There are
in fact still fewer crossings, but we won't need this fact.) In the ($\alpha$) case for $(3)$,
the ($\beta$) and ($\gamma$) cases for (6), and the ($\alpha$) and ($\gamma$) cases for $(7)$, 
the double has a full twist. Still, as before, we find that $M(K_3)\leq 2c(K)-5$, $M(K_6)\leq 2c(K)-5$,
and $M(K_7)\leq 2c(K)-5$. 


Finally, at $(8)$ we have a link, a double of our link $L$ with
a full twist, to which to apply our inductive hypothesis.  By hypothesis, $M(K_8)=M(D(L))=2c(D^\prime) = 2c(L)=2c(K)-2$. 
Having estimates for the $z$-degrees of
the HOMFLY polynomials at all of the leaves of our skein tree diagram, we proceed to work our
way up to the root of the tree, using the inequalities established at the end of section 1. 
Recalling that a vertical arrow represents changing a crossing (passing between $D_-$ and $D_+$)
and a horizontal arrow represents breaking the crossing (passing to $D_0$), and using the labeling
scheme of Figure 11, we find:

\medskip

\noindent $M(K_a)\leq\max\{M(K_1),M(K_A)+1\}\leq\max\{2c(K)-5,2c(K)-5\}=2c(K)-5$ 

\smallskip

\noindent $M(K_b)\leq\max\{M(K_2),M(K_3)+1\}\leq\max\{2c(K)-4,2c(K)-4\}=2c(K)-4$ 

\smallskip

\noindent $M(K_c)\leq\max\{M(K_a),M(K_b)+1\}\leq\max\{2c(K)-5,2c(K)-3\}=2c(K)-3$ 

\smallskip

\noindent $M(K_d)\leq\max\{M(K_4),M(K_5)+1\}\leq\max\{2c(K)-6,2c(K)-4\}=2c(K)-4$ 

\smallskip

\noindent $M(K_e)\leq\max\{M(K_d),M(K_6)+1\}\leq\max\{2c(K)-4,2c(K)-4\}=2c(K)-4$ 

\smallskip

\noindent $M(K_f)\leq\max\{M(K_c),M(K_e)+1\}\leq\max\{2c(K)-3,2c(K)-3\}=2c(K)-3$ 

\smallskip

\noindent $M(K_g)\leq\max\{M(K_B),M(K_7)+1\}\leq\max\{2c(K)-4,2c(K)-4\}=2c(K)-4$ 

\medskip

So far we have used the inequalities; now we establish equalities. 
$M(K_h)\leq\max\{M(K_g),M(K_8)+1\}\leq\max\{2c(K)-4,2c(K)-2\}=2c(K)-2$, but since
$M(K_8)+1=2c(K)-2$ by hypothesis, and $M(K_g)\leq 2c(K)-4$, the two quantities in 
the maximum are unequal, so
$M(K_h)=2c(K)-2$. Finally, $M(D(K))\leq\max\{M(K_f),M(K_h)+1\}\leq\max\{2c(K)-3,2c(K)-1\}=2c(K)-1$,
but since  $M(K_h)+1=2c(K)-1$ and $M(K_f)\leq 2c(K)-3$, the quantities in 
the maximum are again unequal, so $M(D(K))=2c(K)-1$, 
as desired. 

This completes the inductive step, and so the proposition is 
proved by induction, except for dealing with the possibility that the 
links at the leaves of the skein
tree contain twisted doubles of unknots. This is really only an issue at the leaves where
we needed to isotope a short arc through the exterior of the local diagram of the link
before using Morton's inequality to estimate the
$z$-degree, namely at the nodes labeled ``A'',(3),(6),``B'', and (7). In the other cases
we were applying either a local isotopy or our inductive hypothesis to make our
estimates. It is a non-local isotopy, pushing a short arc through the diagram, which 
can expose unlinked components, which would render our
estimates based on the number of crossings of the link inaccurate. However, in every case, 
the link we have after isotopy {\it is} a (possibly twisted) double of a link $L^\prime$, 
obtained, essentially,
by erasing one of the two arcs in the exterior of, but incident to, the original
local diagram, as well as erasing the full twist in the local diagram, replacing the full twist 
with an unknotted arc joining the two remaining ends. Since by hypothesis the
link $L$ is non-split, every unknotted component of $L^\prime$
split off from the rest of $L^\prime$,
which contributes 2 to $M(D(L^\prime))$ without contributing anything to $C(L^\prime)$
comes at the expense of $L^\prime$ having at least two fewer crossings than $L$; the
isotopy of the short arc must ``release'' these unknotted components from the main body 
of $L$. Put differently,
every spanning disk $\Delta$ for the unknotted component must intersect $L$
in its interior, since otherwise the
boundary of a neighborhood of $\Delta$ would provide a splitting sphere for $L$. So the
disk bounding the unknotted, unlinked component, which exists after isotoping the arc, 
must intersect the arc, implying that the
removal of the arc erases at least two crossings in the projection of $L$ (the number
must be even). Consequently, the existence of an unknotted split component, rather 
then requiring us to add 2 to our $z$-degree, implies that it is actually {\it lower} by
2 than our local count might have led us to believe. So in addition to the two crossings
which were lost in the local diagram, which dropped our estimate of the $z$-degree by 4, 
every unknotted split component also, on balance, implies a drop of at least 2 to this
estimate. This implies, at all of the five leaves where this argument is relevant, that
if an unknotted split component is encountered, we can still conclude that 
$M(D(L^\prime))\leq 2c(K)-7$. A quick check above will show that this is as low as
any of the estimates that were used in the arguments above, and so, if anything,
will result in {\it lower} upper bounds than we have used. So it will not affect the final outcome.

\smallskip

With this last detail in place, Proposition 4 is proved.


\heading{\S 3 \\ Non-alternating pretzels do not appear to satisfy deg$_z P_{W(K)}(v,z)=2c(K)$}\endheading

We have shown that for the pretzel knots $K=P(n_1,\ldots ,n_k)$ whose pretzel representation is 
alternating, that is, for which $n_1,\ldots ,n_k\geq 1$, we have $M(W(K))=2g_c(K)$, so $g_c(K)=c(K)=n_1+\cdots +n_k$,
essentially by induction on the $n_i$.
The inductive step, namely the statement of Proposition 4, can be turned around, however, to partially establish
the opposite result for pretzel knots whose pretzel representation is not alternating, that is, the $n_i$ do not all
have the same sign.
In this instance we can show that $M(W(K)) < 2|n_1|+\cdots +2|n_k|$. This follows from the final computation
in the proof of Proposition 4, since if at that point, we know the opposite fact that $M(K_8)$, namely the $z$-degree
for the upper right corner, is {\it less} than $2c(D)-3$, then (since it must be odd) it is at most
$2c(D)-5$, and so, borrowing the notation from the end of the proof, we have 
$M(K_h)\leq\max\{M(K_g),M(K_8)+1\}\leq\max\{2c(D)-4,2c(D)-4\}=2c(D)-4$, so
$M(D(K))\leq\max\{M(K_f),M(K_h)+1\}\leq\max\{2c(D)-3,2c(D)-3\}=2c(D)-3$, so $M(D(K))\leq 2c(D)-3<2c(D)-1$. 
Put differently,
if $M(D(K))=2c(D)-1$, then we {\it must} have $M(K_8)=2c(D)-3$. But arguing our induction in reverse, if we
start with a pretzel knot/link $P(n_1,\ldots ,n_k)$ having (WOLOG) $n_1>0$, $n_2<0$, we can, by turning full
twists into half-twists, return to a base case of $K^\prime=P(1,-1,n_3,\ldots ,n_k)$ to which we can apply a type 2
Reidemeister move, $P(1,-1,n_3,\ldots ,n_k)=P(n_3,\ldots ,n_k)$. So, by the basic construction,
$M(W(P(1,-1,n_3,\ldots ,n_k)))\leq 2|n_3|+\cdots +2|n_k| < 2|1|+2|-1|+2|n_3|+\cdots +2|n_k|= 2c(D^\prime)$, 
where $D^\prime$ is the diagram for $K^\prime$
before the Reidemeister move. Putting back all of the full twists to return to $K$, the above
argument then shows that $M(D(K))< 2c(D)-1$, so $M(W(K)) < 2|n_1|+\cdots +2|n_k|$. Unfortunately,
we are not aware of a proof that $c(P(n_1,\ldots ,n_k))=|n_1|+\cdots +|n_k|$ in the case where this stands
a chance of being true, namely when $\{-1,1\}\notin\{n_1,\ldots ,n_k\}$ .
However, for small examples we do know that this sum is the crossing number.  For example
$P(3,3,-2)$ is the knot $8_{19}$, and $P(3,-3,2)$ is the knot $8_{20}$, and in these cases 
our argument does establish that $M(W(K))<2c(K)$. 

These two examples could of course be computed
directly, and in fact they were, by the authors. It was, in fact, an attempt to show that
$g_c(W(K)) < c(K)$ for these two examples that led the authors to the results presented here. 
Our work here does not establish this, of course; $M(W(K))<2c(K)$ implies that at {\it least}
one of $M(W(K))<2g_c(W(K))$ or $g_c(W(K))<c(K)$ (or both) hold, but to date we do not know which 
of the two might be true! The point to the line of reasoning of the previous paragraph, 
from our perspective, is that
it ``explains'' the inequality $M(W(K))<2c(K)$; it shows that there is a skein tree diagram which
demonstrates that the $z$-degree is ``too low'', without resorting to the need for a fortuitous 
cancellation of terms to occur, to make the degree $2c(K)$ term vanish.

\medskip

On a more positive (or negative) note, the argument above does establish that 
$M(W(K))<2c(K)$ holds for many alternating knots, namely, those that can be formed as a 
connected sum $K_1\# K_2$ for some $K_1\in{\Cal K}$, our class of knots from Proposition 3.
Since these knots are alternating, we know that $c(K_1\# K_2)=c(K_1)+c(K_2)$. 
In this case if we focus on undoing the full twists in the $K_1$ summand, if 
$M(W(K_1\# K_2))=2c(K_1\# K_2)= 2c(K_1)+2c(K_2)$ the argument above will 
allow us to work our way back to $M(W(T_{2,3}\# K_2)) = 2c(T_{2,3})+2c(K_2) =2c(K_2)+6$ by undoing
the full twists introduced to build $K_1$. However, there are still two more full twists in the $T_{2,3}$
summand to work with, and undoing them yields a knot $K$ whose diagram is the projection of $K_2$ with a
single, nugatory, crossing added. Therefore $c(K)=c(K_2)$, but having gotten here in two steps
from $T_{2,3}\# K_2$ we {\it should} have $M(W(K))=2c(K_2)+2$, not $M(W(K))=2c(K_2)$. This contradiction implies
that our original hypothesis, $M(W(K_1\# K_2))=2c(K_1\# K_2)$ is false. Since we must have
$M(W(K_1\# K_2))\leq 2c(K_1\# K_2)$, we must conclude that $M(W(K_1\# K_2)) < 2c(K_1\# K_2)$.
Note that this does {\it not} tell us that $g_c(W(K_1\# K_2)) \neq c(K_1\# K_2)$!
But it does establish that the method of proof used by Tripp, Nakamura, and ourselves to 
show that $g_c(W(K)=c(K)$, namely
to show that $M(W(K))=2c(K)$, will not succeed for all alternating knots.

\medskip

\leavevmode

\epsfxsize=5in
\centerline{{\epsfbox{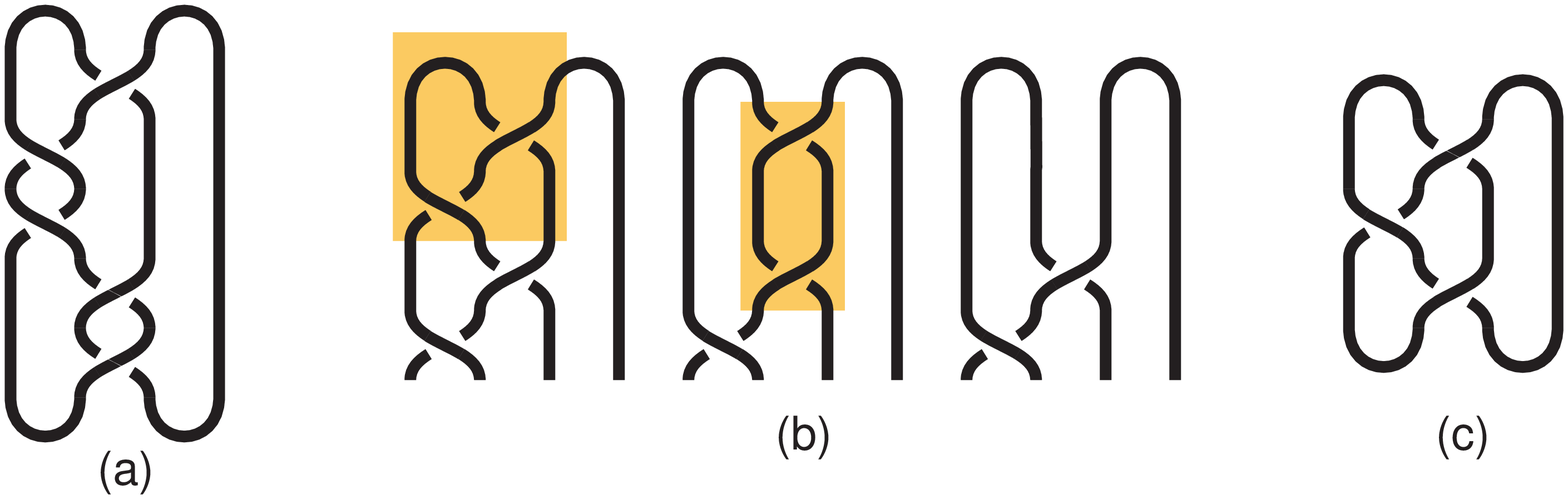}}}

\centerline{Figure 12: Obtaining 2-bridge knots from the trefoil knot}

\medskip

For our final result we establish the last assertion made in the introduction, namely that the class
$\Cal K$ built from the trefoil knots by replacing crossings by full twists already contains the
2-bridge knots, so that including them as initial objects would not enlarge the class. (${\Cal K}$
clearly already contains the $(2,n)$-torus knots.) To see this, we will work backwards from 
an alternating projection of a given 2-bridge knot $K$, successively replacing full twists
with half-twists, and show that we can end up at the trefoil. Reversing this process then 
establishes the result.

The alternating projection for $K$ that we will use will be the alternating 4-plat
projection of [BZ, Proposition 12.13], denoted $\sigma_2^{a_1}\sigma_1^{-a_2}\cdots\sigma_2^{a_m}$ with
$a_i>0$ for $i=1,\ldots ,m$, and illustrated in Figure 12(a). Note that $m$ is necessarily odd. 
The idea is that by 
replacing full twists with half-twists within each block of twists $\sigma_\epsilon^{a_i}$, 
we can steadily lower each of the $a_i$, and
eventually arrive at the knot represented by $\sigma_2^{1}\sigma_1^{-1}\cdots\sigma_2^{1}$. 
Then as shown in Figure
12(b), a further two replacements can, in effect, remove the first pair of $\sigma_2^1\sigma_1^{-1}$'s.
Repeating this process will allow us to reach $\sigma_2^{1}\sigma_1^{-1}\sigma_2^{1}$, shown in
Figure 12(c), which is the (left-handed) trefoil knot. Therefore, reversing this process, 
every 2-bridge knot can be obtained 
from the trefoil knot by a sequence of replacements of crossings by full twists, as desired.
So all 2-bridge knots lie in our class $\Cal K$, and so any knot obtained by replacements 
on a 2-bridge knot will
also lie in $\Cal K$.

\heading{\S 4 \\ Further considerations}\endheading

Proving a result like ours for a class of knots naturally begs the question
``For what other (classes of) knots does $g_c(W(K))=c(K)$ hold?''. The discussion 
of the previous section shows that the techniques we employ cannot establish
this for all alternating knots, since the equality maxdeg$_z(P_{W(K)}(v,z))=2c(K)$ 
need not hold for non-prime alternating knots. It does not, however, show that 
$g_c(W(K))=c(K)$ does not hold.
But based on the results that we have obtained here we feel that it is reasonable to make the

\proclaim{Conjecture}
If $K$ is a nontrivial prime alternating knot, and $W(K)$ is a Whitehead double of 
$K$, then maxdeg$_z(P_{W(K)}(v,z))=2c(K)$, and therefore $g_c(W(K))=c(K)$ .
\endproclaim


The argument given here implies that in trying to establish this conjecture, we may always
inductively replace a pair of crossings forming a bigon by a single crossing. By Euler 
characteristic considerations, the projection of 
any knot, thought of as a graph in the 2-sphere, has a complementary region which is 
either bigon or a triangle. One may think of our results as therefore saying that a proof 
of the conjecture reduces to the cases where the projection has no bigons, i.e., has
complementary regions that all have three or more vertices. This is because our
process of ``deflation'', replacing a full twist with a half-twist, and the inverse process,
keep us within the class of prime, reduced, alternating projections. That no nugatory crossings
can be introduced has been established in Section 2 above. 
Menasco [Me] has shown that
a reduced alternating diagram of a non-prime alternating knot can be detected from
the diagram, that is, there must be a circle in the plane of the projection that 
intersects the projection twice. But such a circle will exist after deflation iff 
there is such a circle before.
If such a circle (before deflation) were to intersect the bigon region, deforming it to 
instead pass through
one of the crossings identifies a nugatory crossing, a contradiction, so the same circle
will work after deflation (showing $\Leftarrow$), while a circle after deflation stays
away from the crossings of the projection, so since the opposite operation of inflation can 
be viewed as local to the crossings (by introducing a very ``small'' bigon), the same circle
also works before deflation (showing $\Rightarrow$). 

\medskip

One may also build
further infinite families of examples, along the lines of Proposition 3, by replacing
the trefoil knot with any other link $K$ for which the equality maxdeg$_z(P_{D(K)}(v,z))=2c(K)-1$
has been established, for example by direct computation. One can, for example, verify
directly that the Borromean rings $L$ satisfies 
maxdeg$_z(P_{D(L)}(v,z) = 11 = 2c(L)-1$, giving rise, using Proposition 4, to a different family of 
(alternating) knots $K$ satisfying $g_c(W(K))=c(K)$ than
the family given by Proposition 3.

Looking beyond the alternating knots, one might look for non-alternating knots $K$ 
to which the techniques of this paper apply, i.e.,
maxdeg$_z(P_{W(K)}(v,z)=2c(K)$.  Our results would then imply that such knots
form the basis for another infinite collection of examples. 
In a different direction, in the course of our proof we find that all of the 
Whitehead doubles of a knot $K$, when our techniques apply to it, have the same canonical genus.
Is there a non-trivial knot $K$ having Whitehead doubles which have different canonical genera?
(The trivial knot does.)
Of course, their HOMFLY polynomials will have the same $z$-degree, so
arguments like ours will be of no help in finding an example.

\medskip

{\bf Acknowledgements:} The first author wishes to thank the Department of Mathematics
of the City College of New York for their hospitality while a part of this work was 
carried out. The second author wishes to thank the Department of Mathematics of the
University of Nebraska - Lincoln and the Nebraska IMMERSE program for their 
hospitality and support.


\enddocument